\documentclass[11pt]{article}
\usepackage{amsmath,amssymb}
\usepackage{theorem}
\usepackage{fullpage}
\usepackage{amsfonts}
\usepackage{latexsym}
\usepackage[english]{babel}
\usepackage{enumerate}
\scrollmode
\newtheorem{corollary}{Corollary}[section]
\newtheorem{theorem}[corollary]{Theorem}
\newtheorem{proposition}[corollary]{Proposition}
\newtheorem{lemma}[corollary]{Lemma}
\theorembodyfont{\upshape}

\newtheorem{definition}[corollary]{Definition}
\newtheorem{remark}[corollary]{Remark}
\newtheorem{example}[corollary]{Example}

\newtheorem{observation}[corollary]{Observation}

\newcounter{abccnt}
\newenvironment{abc}{\begin{list}{\bf(\alph{abccnt})}{\usecounter{abccnt}
\labelwidth4ex \labelsep1ex \leftmargin6ex
\parsep3pt \itemsep1pt \topsep3pt}}{\end{list}}

\newcommand{\barb}{{\rm Barb}}
\newcommand{\ca}{\thickapprox}
\newcommand{\dist}{\not\ca}

\newcommand{\R}{\mathbb{R}}

\newcommand{\F}{\mathbb{F}}
\newcommand{\Z}{\mathbb{Z}}

\newcommand{\rz}{\mathbb{R}}
\newcommand{\cz}{\mathbb{C}}
\newcommand{\gz}{\mathbb{Z}}

\newcommand{\abs}[1]{\lvert#1\rvert}

\newcommand{\whom}{\mathrm{w}_{\mathrm{hom}}}

\DeclareMathOperator{\GF}{GF}

\DeclareMathOperator{\OA}{OA}

\DeclareMathOperator{\Hom}{Hom}
\DeclareMathOperator{\Aut}{Aut}
\DeclareMathOperator{\soc}{Soc}
\DeclareMathOperator{\rad}{Rad}
\DeclareMathOperator{\kernel}{Ker}

\newcommand{\qed}{\nolinebreak\hfill$\Box$\vspace{3pt}}
\newcommand{\pf}{{\em Proof\/: }}

\title{Ring Geometries, Two-Weight Codes,\\ and Strongly Regular Graphs}
\author{Eimear Byrne, Marcus Greferath, and Thomas Honold}

\date{}

\begin{document}
\maketitle

\begin{abstract}
  It is known that a linear two-weight code $C$ over a finite
  field $\F_q$ corresponds both to a multiset in a projective space over
  $\F_q$ that meets every hyperplane in either $a$ or $b$ points for
  some integers $a<b$, and to a strongly regular graph whose vertices
  may be identified with the codewords of $C$.
 
  Here we extend this classical result to the case of a ring-linear code
  with exactly two nonzero homogeneous weights and multisets of points in    
  an associated projective ring geometry. We will
  show that a two-weight code over a finite Frobenius ring gives rise to a
  strongly regular graph, and we will give some constructions of two-weight
  codes using ring geometries. These examples all yield infinite families of 
  strongly regular graphs with non-trivial parameters.
\end{abstract}

\section{Introduction}

  It has been known since the 1970's that there are connections
  between certain sets in a projective geometry, partial difference
  sets, codes with two nonzero weights, and strongly regular
  graphs. The equivalence of classes
  of these objects was first observed by Delsarte \cite{D71,D72}, and since
  then, many other constructions have been discovered (see
  \cite{CK86} for a survey). 

  More recently, it has been observed that a linear code over $\Z_4$
  with exactly two nonzero Lee weights also determines a strongly
  regular graph \cite{BF2010}, so that a correspondence may be drawn
  between certain distance invariant non-linear binary codes and such
  graphs. In this paper, we show that a linear code over a finite
  Frobenius ring with two nonzero homogeneous weights also determines
  a strongly regular graph. In addition, we give constructions of
  two-weight codes over finite chain rings determined by point sets of
  associated projective geometries. These point sets provide new
  constructions of non-trivial strongly regular graphs. This indicates that 
  projective geometries over rings may be used to
  determine new families of strongly regular graphs.

  \section{Preliminaries: Finite Frobenius rings and homogeneous
    weights}

  The merits of studying properties of linear codes over finite rings has 
  been recognized since
  the end of the eighties (cf.~\cite{nech91,H94}). Two important
  foundational results by J.\ Wood \cite{W99,Wa99} state that
  the MacWilliams' equivalence and duality theorems 
  can be extended for linear codes over a very large class of
  rings, namely the finite Frobenius rings. 
  A recent paper by Wood~\cite{wood2007} yields the final justification for the belief that the finite Frobenius rings are the appropriate class for algebraic coding theory over rings. The reader should note however, that the Frobenius property of the ring can be dropped by passing to a suitable module (a quasi-Frobenius module of the given ring) serving as alphabet for ring-linear coding theory (cf.~\cite{grefnechwisb02}).    
  
 Detailed from the theory of quasi-Frobenius and Frobenius rings can be found in \cite{F67,F76,lam}. Here, we briefly describe properties of such rings in
  the finite case.

  For a finite ring $R$ denote by $\widehat{R}:=\Hom_\gz(R,\cz^\times)$ the character module of $R$.   
  $\widehat{R}$ has a natural $R$-bimodule structure determined by
  $\chi^r(x) := \chi(rx)$ and ${^r\chi}(x) := \chi(xr)$ for all $r,x\in
  R$ and each character $\chi \in \widehat{R}$.

  A finite ring $R$ is called a \emph{Frobenius ring} if it satisfies any
  of the following equivalent conditions:
  \begin{itemize}
  \item[(i)]  $_R(\soc({}_RR))\cong{}_R(R/\rad(R))$;
  \item[(ii)] $(\soc(R_R))_R\cong(R/\rad (R))_R$;
  \item[(iii)] ${}_RR\cong{}_R\widehat{R}$;
  \item[(iv)] $R_R\cong\widehat{R}_R$;
  \item[(v)] $\soc({}_RR)$ is a principal left ideal;
  \item[(vi)] $\soc(R_R)$ is a principal right ideal.
  \end{itemize}
  Conditions~(iii) and~(iv) say that $R$ has a
  \emph{left (resp.~right) generating character} $\chi$, which means $\widehat{R}=\{{}^r\chi\mid r\in R\}$ (resp.~$\widehat{R}=\{\chi^r\mid r\in R\}$). It can be shown that a character is left generating if and only if it is right generating.

  Homogeneous weights were first introduced by I.~Constantinescu and
  W.~Heise \cite{consheis97} in the context of linear codes over
  integer residue rings. Such a weight function has the property of
  taking a constant value on sets of associated ring elements.

  Weights with these properties have been revisited and generalized in
  different ways (cf.~\cite{gray,bruce,thomas-nech99,werner-thomas00}). Here, we
  follow the line in \cite{gray,bruce} which works without
  restrictions on the underlying finite ring. 

\begin{definition}\label{defhomogen} Let $R$ be a finite ring and
  let $w: R \Longrightarrow \R$ be a map satisfying $w(0)=0$ and $w(x) \geq 0$ for all $x\in R$.
  The map $w$ is
  called a \emph{(left) homogeneous weight}, if $w(0)=0$ and the following is
  true:
  \begin{abc}
  \item[\bf (H1)] If $Rx=Ry$ then $w(x)= w(y)$ for all $x,y\in R$.
  \item[\bf (H2)] There exists a real number $\gamma$ such that
    \begin{equation*}
  \sum_{y\in Rx}w(y) \; =\; \gamma \, |Rx|\qquad\text{for all $x\in
    R\setminus \{0\}$}.
\end{equation*}
\end{abc}
\end{definition}

The number $\gamma$ may be thought of as the {\em average value}
of the weight function, and the condition {\bf (H2)} simply states
that this average is the same on all nonzero principal left ideals.
 
\begin{example}\label{exhw}
  We give some examples of finite rings with homogeneous weight functions. Observe that in 
  \ref{F2}, the a non-zero ring element has been assigned weight zero. 
  \begin{abc}
  \item On every finite field $\F_q$ the Hamming weight
    is a homogeneous weight of average value $\gamma=\frac{q-1}{q}$.
  \item On $\gz_4$ the Lee weight is homogeneous with
    $\gamma=1$.
  \item On a finite chain ring $R$ with residue field $\GF(q)$ the
    weight
    \begin{equation*}
      w\colon R\longrightarrow \rz,\quad x \mapsto
      \begin{cases}
        0&\text{if $x=0$},\\
        \frac{q}{q-1}&\text{if $x\in\soc(R)$, $x\neq 0$},\\
        1&\text{otherwise},
      \end{cases}
    \end{equation*}
is a left and right homogeneous weight with $\gamma=1$.
\item \label{F2} On the ring $\F_2\oplus\F_2$ the weight $w$ defined by
  $w(0,0)=w(1,1)=0$, $w(1,0)=w(0,1)=2$ is a homogeneous weight with $\gamma=1$.
\item On the (local) ring $R=\F_2[X,Y]/(X^2,Y^2, XY)$ 
  of order $8$ the weight
  \begin{equation*}
    w\colon R\longrightarrow \rz,\quad x \mapsto
  \begin{cases}
    0&\text{if $x=0$},\\
        2&\text{if $x\in\rad(R)\setminus\{0\}=\{X,Y,X+Y\}$},\\
        \frac{1}{2}&\text{otherwise},
  \end{cases}
\end{equation*}
is a homogeneous weight with $\gamma=1$.
\end{abc}
\end{example}

As is common in coding theory, a weight $w$ on a finite ring $R$ will be extended
to a weight on the $R$-module ${}_RR^n$ in the obvious way:
\begin{equation*}
w\colon R^n \longrightarrow \R,\quad w(c_1,\dots,c_n) = \sum_{i=1}^n w(c_i).  
\end{equation*}

It can be shown that up to the choice of $\gamma$, every finite ring admits a unique (left) homogeneous weight (cf.~\cite[Thm.~1.3]{bruce}). For the investigations here a further observation by Honold \cite{thomas01} is important and is used in the proof of Proposition \ref{prop3}.

\begin{proposition}\label{char-homogen} 
  Let $R$ be a finite Frobenius
  ring with generating character $\chi$.
  Then the (left) homogeneous weight on $R$ has the form
  \begin{equation*}
w\colon R\longrightarrow \rz, \quad x \mapsto
  \gamma\Big[1-\frac{1}{|R^{\times}|}\sum_{u\in R^{\times}}
  \chi(ux)\Big]
\end{equation*}
for some real number $\gamma$.
\end{proposition}

The following important fact, stated in \cite[Lem.~1.5,
Rem.~1.7(b)]{bruce} and will be used for our computations later (see Lemma \ref{colweight}).

\begin{remark}\label{mufrob} 
  If $w$ is a (left) homogeneous weight of
  average value $\gamma\neq 0$ on the finite ring $R$, then $R$ is a
  Frobenius ring if and only if:
  \begin{abc}
  \item[\bf(H2*)] 
    \begin{equation*}
  \sum_{y\in I}w(y) \; =\; \gamma \, |I|\qquad\text{for all $\{0\}
  \neq I \leq {_RR}$}.
\end{equation*}
  \end{abc}
\end{remark}

Note, for example, that the ring of Ex.~\ref{exhw}\,(e) is not a Frobenius ring, since its socle is not cyclic.
Indeed, $$\sum_{x\in\rad(R)}w(x)=6\neq\abs{\rad (R)}$$ for the given weight.

It has been shown in ~\cite{thomas-nech99} that if $R$ is a finite
Frobenius ring, then every left homogeneous weight is also right
homogeneous, clearly with the same average value $\gamma$.
This can also be inferred from Proposition~\ref{char-homogen} and
the equation $$\sum_{u\in
  R^{\times}}\chi(ux)\; = \; \sum_{{\psi\in\widehat{R} \atop \psi\mbox{ \footnotesize generating }}}
\psi(x)\; = \; \sum_{u\in
  R^{\times}}\chi(xu).$$

We shall call the unique homogeneous weight of average value
$\gamma=1$ the {\em normalized\/} homogeneous weight and denote it by
$\whom$.

Finally, the proof of Lemma \ref{lem2} uses the notion of  \emph{M\"obius inversion} 
for functions on finite posets\footnote{A poset is a set together with a partial order.} 
For details see \cite[Ch.~IV]{aign97},~\cite{rota64},~\cite[Ch.~3.6]{stan97}, or \cite{wieg59}.

Let $P$ be a finite poset. Consider the function $\mu:
P\times P \longrightarrow \rz$, defined recursively by the equations
$\mu(x,x) = 1$ and
\begin{equation*}
\mu(y,x) =
\begin{cases}
    - \sum\limits_{y < t \leq x} \mu(t,x)&\text{if $y<x$},\\
    0 &\text{if $y\not\leq x$}.
\end{cases}
\end{equation*}

The map $\mu$ is called the \emph{M\"obius function} of $P$ and induces
for arbitrary pairs of real-valued functions $f,g$ on $P$ the
equivalence:
\begin{equation*}
g(x) = \sum_{y\leq x}f(y)\quad\text{for all $x\in P$}\qquad\iff\qquad
f(x) = \sum_{y\leq x}g(y)\,\mu(y,x)\quad\text{for all $x\in P$}.
\end{equation*}

This equivalence is usually referred to as {\em M\"obius inversion\/}. 

For any finite ring $R$ we denote by $R^\times$ the group of units of $R$.
An application of M\"obius inversion pertinent to this paper is then given by the following.

\begin{example}\label{moeb}
  Let $R$ be a finite ring and let $\mu$ denote the M\"obius function on the poset 
  $\{Rx\mid x\in R\}$, partially ordered by set inclusion. 
  Let $x\in R$. Then $R^\times x$ is the set of all generating elements of
  of $Rx$ and
  \begin{equation*}
    |R^\times x|=\sum_{Ry\leq Rx}|Ry|\,\mu(Ry,Rx)
 \end{equation*}
 holds since $\displaystyle{|Rx| = \sum_{Ry \leq Rx} |R^\times y |}$.
\end{example}

\section{Linear codes and geometries over finite Frobenius rings}

Let $R$ be a finite ring $R$.  A left linear code of length $n$ over $R$ is
a submodule of ${}_RR^n$, which we indicate by writing $C\leq {_RR^n}$. We say that
$C$ is {\em $k$-generated\/} if $C$ possesses a generating set of $k$ elements. Equivalently, $C$ is the row space of a $k\times n$-matrix over $R$. Note that a $k$-generated code is always $\ell$-generated for all $\ell \geq k$.

A pair of left linear codes $C$ and $C'$ over $R$ are said to be \emph{isometric} relative to the weight function $w$ if there exists a bijective linear $w$-isometry from $C$ onto $C'$, i.e.~an $R$-linear isomorphism between the codes that preserves $w$. If $w$ is homogeneous and $R$ is a finite Frobenius ring this occurs precisely when the given isomorphism can be extended to a monomial transformation from ${}_RR^n$ onto itself (cf.~\cite{bruce}). 

For what follows, we need to introduce the notion of projective geometry over a ring. There are various approaches to this in the literature. Here we follow Veldkamp~\cite[26.4.1]{veld} and adopt the notion of a {\em Barbilian space\/} of a right module $M_R$. Let $M^{*}:= {\rm Hom}({M_R},{R_R})$ denote the dual of $M_R$. $M^{*}$ forms a left $R$-module by the right action of $R$ on $M$. 

The {\em Barbilian space of
 $\barb(M_R)$} is defined as the quadruple $(P,H,\mid,\dist)$ where
\begin{eqnarray*}
   P & := & \{ xR \mid x\in M {\rm\ with\ } \varphi x = 1 {\rm\ for\ some\ }
             \varphi \in M^{*}\}\\
   H & := & \{ R\varphi \mid \varphi \in M^{*} {\rm\ with\ } \varphi x = 1
             {\rm\ for\ some\ } x \in M\}\\
   | & := & \{ (R\varphi,xR) \mid (\varphi,x) \in M^*\times M {\rm\ with\ }
             \varphi x = 0\} \, \cap \, H \times P\\
   \dist & := & \{ (R\varphi,xR) \mid (\varphi,x) \in M^*\times M {\rm\ with\ }
             \varphi x = 1\}
\end{eqnarray*}

The elements of $P$ are called {\em points\/}, the elements of $H$ are called {\em hyperplanes\/} of $\barb(M_R)$. Relation $|$ is called {\em incidence\/}, and relation $\dist$ is called {\em distant\/}. 

A left linear code $C$ may now be viewed as a set of evaluations of an
$n$-tuple of linear functionals $\psi_i: {_RR}^k \longrightarrow {}_RR$,
$i=1,\dots,n$, for some positive integer $k$\,:
\begin{equation*}
C = \left\{\bigl(x\psi_1,x\psi_2,...,x\psi_n\bigr)\mid x \in R^k\right\}.
\end{equation*}
A generator matrix for $C$ is then given as the $k \times n$ generator matrix
$Y=(y_1|y_2|\dots|y_n)$, with each $y_i \in R^k$, if the linear
functionals $\varphi_i$ correspond to taking the standard inner product
$x\psi_i := x\cdot y_i $ for each $x \in R^k$. Let $\pi_i: R^n \longrightarrow R,\;\; x\mapsto x_i$ denote the projection of $R^n$ onto its $i$-th coordinate. If all coordinate projections $\pi_i$ are surjective, then up to monomial equivalence, the code $C$ is equivalently described by the multiset $\{y_1R,\dots,y_nR\}$ of points in $\barb(R^k_R)$. Such a code will be called a {\em regular\/} code.
 
For the purposes of this paper, we will always require that the projection of $C$ on any coordinate returns the full ring $R$, so that each column of the generator matrix $Y$ generates a point in $\barb(R^k_R)$. We will also usually require that the points $\{y_1R,...,y_nR\}$ are all distinct, in which case we will refer to the generated code as a {\em projective\/} code. 

Recall that there exist finite Frobenius rings such that the homogeneous weight of a nonzero element is zero, in which case a linear code over that ring could have nonzero words of weight zero. In most what follows we will rule out such an anomaly of a code by requiring that the code be {\em proper\/}.

More formally we give the following definition.

\begin{definition}\label{defproj}
  A linear code $C \leq {_RR^n}$ with $k \times n$ generator matrix
  $Y=(y_1|y_2|\dots|y_n)$ is called
  \begin{enumerate}
  \item[(i)] \emph{regular} if $y_iR$ is a point of $\barb(R^k_R)$ for each $i \in
    \{1,\dots,n\}$,
  \item[(ii)] \emph{projective} if $y_iR\neq y_jR$ for any pair of distinct
    coordinates $i,j \in \{1,\dots,n\}$,
   \item[(iii)] \emph{proper} if $\whom(c)>0$ for every nonzero $c \in C$.   
  \end{enumerate}
\end{definition}

\begin{remark}
  \label{1-r-s-rs}
  A linear code $C\leq{}_RR^n$ is projective if and only if
  the $n$ coordinate projections $\pi_i: C\longrightarrow R$, $c\mapsto c_i$ generate
  distinct cyclic submodules of $_R\Hom({}_RC,{}_RR)$. Hence
  projectivity of a linear code is a well-defined concept and independent of the particular choice of $Y$. 
\end{remark}

\section{Two-weight codes over a finite Frobenius ring}

In this section we introduce codes with two non-zero homogeneous weights and
make some basic observations. In what follows let $R$
be a finite Frobenius ring, and let $w$ be a homogeneous weight
on $R$ of average value $\gamma\neq 0$.

The following can be deduced directly from the homogeneity
condition {\bf H2*}. It is an immediate extension of a statement given in
\cite{consheis97} (cf.~also \cite{bruce}).

\begin{lemma} \label{colweight} If $C\leq {_RR^n}$ is a linear code
  over $R$ then for all $i=1, \ldots, n$ there holds
  \begin{equation*}
    \sum_{c\in
    C}w(c_i) \; = \;
  \begin{cases}
  \gamma\,|C|&\text{if $\pi_i(C) \neq 0$},\\
      0&\text{otherwise}.
\end{cases}
\end{equation*}
\end{lemma}

\begin{definition}\label{two-weight}
  A proper linear code $C \leq {_RR^n}$ is called a \emph{two-weight code},
  if $w$ takes exactly two nonzero values $w_1$ and $w_2$ on $C$. 
\end{definition}

If the normalized homogeneous weight $\whom$ on $R$ is strictly positive, i.e.~if $\whom(x)>0$ for all $0\neq x\in R$, then all codes over $R$ are proper. 

The class of finite Frobenius rings where the normalised homogeneous weight is not strictly positive was characterised in \cite[Th.~1]{thomas-nech99}:

  \begin{proposition}\label{R_0}
    The normalized homogeneous weight $\whom$ on a finite Frobenius
    ring $R$ is positive definite if and only if $R$ has at
    most one two-sided ideal of cardinality $2$. 
  \end{proposition}

Since the sum of the weights in a nonzero coordinate of $C$ is fixed and given by $\gamma |C|$, we have the following system of equations.

\begin{observation}\label{linearsystem}
  Let $C$ be a two-weight code of length $n$ with weights $w_1$ and $w_2$. Let $b_1$ denote the number of codewords of weight
  $w_1$ and $b_2$ the number of codewords of weight $w_2$. If $\pi_i(C) \neq \{0\}$ for all $i=1, \ldots, n$ then
  \begin{equation*}
    b_1 w_1 + b_2 w_2 \;  = \;  \gamma\, n\,|C| \quad\text{and}\quad
    b_1 + b_2 \; = \; |C|-1.
\end{equation*}
\end{observation}

This system of equations determines the relationship between the
$b_i$, $|C|$, $w_i$ and $n$. In particular, it shows that the values
$b_i$ are uniquely determined by $n$, $\abs{C}$, $w_1$ and $w_2$.  We
prefer to rewrite the above system in the following more convenient
matrix form:
\begin{equation*}
  \begin{bmatrix}
    w_1 & w_2\\ 1 & 1\end{bmatrix}\cdot
  \begin{bmatrix}
    b_1 \\b_2
  \end{bmatrix}\; = \;
  \begin{bmatrix}
    \gamma\,n\, |C|\\|C|-1
  \end{bmatrix}.
\end{equation*}
 
Given a linear code $C\leq{}_RR^n$ and a left ideal $I$ of $R$, for each $i \in \{1,...,n\}$ we write 
$C(i,I)$ to denote the left $R$-submodule of $C$ defined by
$$C(i,I): = \{ c \in C \:|\: \pi_i(c) \in I\}.$$  

\begin{lemma} 
  \label{lem1} 
  Let $C\leq{}_RR^n$ be a linear code and
  let $I$ be a left ideal of $R$. Then for each $i \in \{1,\dots, n\}$,
  \begin{equation*}
  \abs{C(i,I)} \; = \; \frac{|I|}{|I +
    \pi_i(C)|}\, |C|.    
  \end{equation*}
  In particular, if $\pi_i(C) = R$ then
  \begin{equation*}
  \abs{C(i,I)} \; = \; \frac{|I|}{|R|} \, |C|.    
  \end{equation*}
\end{lemma}

\pf First note that $\abs{C(i,I)} = \sum_{r\in I}
\abs{\{c\in C\mid\pi_i(c) = r\}}$. Now if $r$ is an element in $I \cap
\pi_i(C)$ then $r$ occurs in the $i$-th position of exactly
$\abs{\kernel(\pi_i)\cap C}$ words of $C$. Then, summing over all $r$
in $I$, we obtain
\begin{equation*}
\abs{C(i,I)}\; = \;
\abs{I\cap \pi_i(C)} \cdot \abs{\kernel(\pi_i)\cap C}.
\end{equation*}
A simple homomorphism argument completes the proof.\qed

An immediate consequence of Lemma \ref{lem1} is that if $C$ is
regular, then the expression $|C(i,I)|$ depends
on $C$ and $I$, certainly, but is independent of the chosen coordinate
$i$. 

\begin{lemma}
  \label{c_ic_j}
  Given a regular projective code $C\leq{}_RR^n$ and any pair
  $i,j\in\{1,\dots,n\}$ of distinct coordinates, there exists a
  codeword $c=(c_1,\dots,c_n)\in C$ with $c_i=0$ and $c_j\neq 0$.
\end{lemma}
\pf Let $C$ have $k \times n$ generator matrix $(y_1|\cdots| y_n)$.
Since $C$ is regular, $\pi_i(C)=y_iR$ for each $i$ and $C/\kernel \pi_i = y_iR$. Since 
$C$ is projective the $y_iR$ are distinct and hence 
$C \cap \kernel\:\pi_i \neq C\cap \kernel\:\pi_j$ for $i \neq j$. 
Furthermore, since $C$ is regular, from Lemma \ref{lem1}
$\abs{C\cap \kernel\:\pi_i}=\abs{C \cap \kernel\:\pi_j }=\frac{\abs{C}}{\abs{R}}$. Now, the result follows.\qed

Given a two-weight code $C$ with weights $w_1$ and $w_2$ define
\begin{eqnarray}\label{eqc1c2}
  C_1 \; := \; \{c\in C \mid w(c) = w_1\}, \;\;
  \mbox{and}\;\; C_2 \; := \; \{c\in C\mid w(c)= w_2\}.
\end{eqnarray}
For an arbitrary subset $T\subseteq R$ we define $$b_1 (i,T) \; := \; |\{c\in C_1 \mid
\pi_i(c) \in T\}|,\;\; \mbox{and}\;\; b_2 (i,T) \; := \; |\{c\in C_2
\mid \pi_i(c) \in T\}|.$$

Our next goal is to show that $b_1(i,I)$ (and hence $b_2(i,I)$) does
not depend on the chosen coordinate $i$, provided $C$ is regular and
projective and $I$ is a left ideal of $R$. 

Using Lemma~\ref{lem1} we extend Observation~\ref{linearsystem} in the
case of regular projective codes as follows.

\begin{proposition} \label{4.6} Let $C\leq {_RR^n}$ be a regular projective
   two-weight code with weights $w_1$ and $w_2$. For each $i\in \{1,
  \dots,n\}$ and $I \leq {}_RR$ then
  \begin{equation*}
    \begin{bmatrix}
      w_1 & w_2\\ 1 & 1
    \end{bmatrix}\cdot
    \begin{bmatrix}
      b_1 (i,I)\\b_2 (i,I)
    \end{bmatrix}\; = \;
    \begin{bmatrix}
      \gamma\, n\, \frac{|I|}{|R|}\,
     |C|\\\frac{|I|}{|R|}\, |C|-1
    \end{bmatrix}
\end{equation*}
 if $I \neq 0$,
 and
 \begin{equation*}
   \begin{bmatrix}
     w_1 & w_2\\ 1 & 1
   \end{bmatrix}\cdot
   \begin{bmatrix}
     b_1 (i,0)\\b_2 (i,0)
   \end{bmatrix}\; = \;
   \begin{bmatrix}
     \gamma\,
      (n-1)\,\frac{|C|}{|R|}\\\frac{|C|}{|R|} -1
   \end{bmatrix}
\end{equation*}
  provided $I=0$. In particular, the values $b_1 (i,I)$ and $b_2 (i,I)$ are
  independent of $i$ for all nonzero $I\leq {_RR}$.

\end{proposition}

\pf We count in two ways the total homogeneous weight of the linear
subcode $C(i,I)$ of $C$. From
Lemma~\ref{colweight} we have 
\[w_1b_1(i,I) + w_2b_2(i,I)=\sum_{c\in
  C(i,I)}w(c)=\sum_{j=1}^n\sum_{c\in C(i,I)}w(c_j)=\gamma s|C(i,I)|,\]
where $s:=\left|\left\{j\mid\pi_j\bigl(C(i,I)\bigr)\neq
    0\right\}\right|$. Since $C$ is regular and projective, from Lemma~\ref{c_ic_j} we have
$\pi_j\bigl(C(i,\{0\})\bigr)\neq 0$ for $i\neq j$.
Hence, $s=n-1$ or $s=n$ depending on whether $I=0$ or
$I\neq 0$, respectively. Using Lemma~\ref{lem1} we obtain
$|C(i,I)|=\frac{|I|}{|R|}|C|$ for any $I$, which yields the result in both
cases.\qed

In light of the fact that the numbers $b_1(i,I)$ and $b_2(i,I)$ are independent of the choice of $i$ for a given regular projective two-weight code, we will refer to these magnitudes simply as $b_1(I)$ and $b_2(I)$. Of course, $b_1(R) = b_1 $ and $b_2(R) = b_2 $ as introduced in Observation~\ref{linearsystem}.

\section{Two-weight codes and strongly regular graphs}

We will now draw connections between two-weight codes and strongly
regular graphs.  In particular, we show that every linear code over a
finite Frobenius ring with exactly two non-zero homogeneous weights
renders a strongly regular graph.

\begin{definition}
  A simple graph $\Gamma = (V,E)$ with vertex set $V$ and edge set $E$ is
  called strongly regular with parameters $(N,K,\lambda,\mu)$ if:
  \begin{itemize}\itemsep=0mm
  \item[(i)] $\Gamma$ has $N$ vertices and is neither empty nor complete, 
  \item[(ii)] $\Gamma$ is regular of degree $K$, and 
  \item[(iii)] Every adjacent pair of vertices $v,v' \in V$ has exactly $\lambda$
  common neighbours in $V$
  \item[(iv)] Every non-adjacent pair of vertices $v,v' \in V$ has exactly $\mu$
  common neighbours in $V$.
  \end{itemize}
\end{definition}

Strongly regular graphs are well-studied and are equivalent to
symmetric association schemes with $2$ classes. Corresponding
to any strongly regular graph $\Gamma$ is a 3-dimensional Bose-Mesner
algebra spanned by the $N\times N$ identity matrix $I$, the all-one
matrix $J$ and the adjacency matrix $A$ of $\Gamma$.  The reader is
referred to \cite{CK86,CL91,lw98} for properties of such graphs and
their relations to other combinatorial objects.  The complement of a
strongly regular graph with parameters $(N,K,\lambda,\mu)$ is also
strongly regular and has parameters
$(N,N-K-1,N-2K+\mu-2,N-2K+\lambda)$.  A strongly regular
graph is said to be trivial if either $\Gamma$ or its complement is a
disjoint union of complete graphs (of the same size). A strongly regular graph
$\Gamma$ is called nontrivial if and only if its parameters satisfy the
condition
 $ 0 < \mu < K$.

The parameters of a strongly regular graph must satisfy
certain {\em feasibility conditions\/}. For example, if $\Gamma = (V,E)$
is strongly regular with parameters $(N,K,\lambda,\mu)$ then a simple
counting argument shows that $K(K-\lambda-1) = \mu(N-K-1)$. The parameters
$N,K,\lambda,\mu$ are called feasible if they satisfy well-known feasibility conditions, and an important aspect of the theory of strongly regular graphs is to establish the existence of a strongly regular graph for a given set of feasible parameters.

Let $G$ be an additive abelian group and let $D$ be a subset of $G$
such that $-d \in D$ for each $d \in D$ and $0 \notin D$. The {\em Cayley graph} 
of $G$ with respect to $D$ is given by the graph
$\Gamma=(G,E)$ with vertex set $G$ and edge set $E
=\bigl\{(g,d+g)\mid g \in G, d \in D\bigr\}$. 
Then $\Gamma$ is regular of degree $|D|$ and
is strongly regular if and only if $D$ is a so-called (regular)
partial difference set (see ~\cite{CK86,ma94}).

In this section we consider this construction for the particular case where $G$
is a submodule of ${}_RR^n$ (i.e. an $R$-linear code) and $D$
is a subset of $G$ (in fact the set of words of a constant weight subcode of $G$).
We will show that if $G$ is a regular projective
two-weight code, and $D$ consists of the codewords of a given weight, say
$w_1$, then the graph $\Gamma=(G,E)$ defined above is strongly regular.

\begin{lemma}\label{lem2}
  Let $C\leq {}_RR^n$ be a regular projective two-weight code with
  weights $w_1$ and $w_2$, and let $C_1$ and $C_2$ be as defined in
  Equation~\eqref{eqc1c2}.  Then for each $i \in \{1,\dots,n\}$ we
  have
  \begin{equation*}
\sum_{x\in C_1}w(x_i) \; =\; \frac{b_1\, w_1}{n}\qquad
\text{and}\qquad\sum_{x\in C_2}w(x_i) \; =\; \frac{b_2\, w_2}{n}.
\end{equation*}
\end{lemma}

\pf We need only show $\sum_{x\in C_1}w(x_i)$ is independent of
$i$, since in that case we find that the total homogeneous weight of
the words of weight $w_1$ in $C$ (which is given by $w_1b_1$) equals
$n\sum_{x\in C_1}w(x_i)$.  As in Example \ref{moeb}, first observe that  
$$b_1(i, R^\times x) \; = \; \sum_{Ry\leq Rx}
b_1(i, Ry) \, \mu(Ry,Rx)$$ by M\"obius inversion. By Proposition \ref{4.6}, these magnitudes clearly do not depend on the chosen coordinate $i$, and so 
we can write $b_1(i,R^\times x) = b_1(R^\times x)$ for all $i$. Finally we obtain \begin{equation*}
\sum_{x\in C_1}w(x_i) \; = \; \sum_{Ry\leq R}
w(y)\, b_1(R^\times y)   
\end{equation*}
which again will not depend on the choice of $i$. This completes the proof for both claims.\qed

\begin{proposition}\label{prop3}
  Let $C\leq {_RR^n}$ be a regular projective two-weight code with
  weights $w_1$ and $w_2$. Then
  \begin{equation*}
    \sum_{x\in C_1} w(x-c) \; = \;\gamma\,n\,b_1\Big[1 -
    \Big(1-\frac{w_1}{\gamma\, n}\Big)\Big(1-\frac{w(c)}{\gamma\,
      n}\Big)\Big]    
  \end{equation*}
for every $c\in C$.
\end{proposition}

\pf Using Proposition \ref{char-homogen} we compute $\sum_{x\in C_1}
w(x-c)$ as follows.
\begin{align*}
  \sum_{x\in C_1} w(x-c) & = \sum_{x\in C_1}\sum_{i=1}^n w(x_i-c_i)\\
  & = \gamma\, \sum_{x\in C_1}\sum_{i=1}^n
  \Big[1-\frac{1}{|R^\times|}\sum_{u\in R^\times}\chi(u[x_i-c_i])\Big]\\
  & = \gamma \,n\, b_1 -\gamma\sum_{x\in
    C_1}\sum_{i=1}^n\frac{1}{|R^\times|}\sum_{u\in
    R^\times}\chi(-uc_i)\chi(ux_i)\\
  & = \gamma \, n\,b_1
  -\gamma\sum_{i=1}^n\frac{1}{|R^\times|}\sum_{u\in
    R^\times}\chi(-uc_i)\Big[\sum_{x\in C_1}\chi(ux_i)\Big].
\end{align*} 
Observing that as $uC_1=C_1$ for all $u\in R^\times$, we get
\begin{equation*}
\sum\limits_{x\in
  C_1}\chi(ux_i) \; = \; \sum_{x\in C_1}\chi(x_i)\; = \;
\frac{1}{|R^\times|}\sum_{v\in R^\times}\sum_{x\in
  C_1}\chi(vx_i)\; = \; \sum_{x\in
  C_1}\frac{1}{|R^\times|}\sum_{v\in R^\times}\chi(vx_i), 
\end{equation*}
and so we may write
\begin{align*}
  \sum_{x\in C_1} w(x-c) & = \gamma \, n\, b_1
  -\gamma\sum_{i=1}^n\Big[\frac{1}{|R^\times|}\sum_{u\in
    R^\times}\chi(-uc_i)\Big]\sum_{x\in C_1}
  \Big[\frac{1}{|R^\times|}\sum_{v\in R^\times}\chi(vx_i)\Big]\\
  & = \gamma \,n\, b_1
  -\gamma\sum_{i=1}^n\Big(1-\frac{w(c_i)}{\gamma}\Big)\sum_{x\in C_1}
  \Big(1-\frac{w(x_i)}{\gamma}\Big).
\end{align*}
Using Lemma~\ref{lem2} we can rewrite
the latter expression as $$\sum_{x\in C_1}
\Big(1-\frac{w(x_i)}{\gamma}\Big) \; = \; b_1 - \frac{b_1\,
  w_1}{\gamma\,n} \; = \; b_1\Big(1- \frac{w_1}{\gamma\,n}\Big),$$
which is clearly independent of the coordinate $i$, so that we have
\begin{align*}
  \sum_{x\in C_1} w(x-c) & = \gamma \,n\, b_1 -\gamma\, b_1 \Big(1
  - \frac{w_1}{\gamma\,n}\Big)\sum_{i=1}^n\Big(1-\frac{w(c_i)}{\gamma}\Big)\\
  & = \gamma \,n\,  b_1 - \gamma\,n\, b_1 \,\Big(1-\frac{w_1}{\gamma\,
    n}\Big)\Big(1 - \frac{w(c)}{\gamma\, n}\Big)\;\\
  & = \gamma\,n\, b_1
  \Big[1 - \Big(1-\frac{w_1}{\gamma\, n}\Big)\Big(1-\frac{w(c)}{\gamma\,
    n}\Big)\Big].
\end{align*}
This completes the proof.\qed

This result allows us to set up a new system of equations, as
indicated below.

\begin{corollary}\label{corde}
  Let $C\leq {_RR^n}$ be a regular projective two-weight code with
  weights $w_1$ and $w_2$, and let $C_1$ and $C_2$ be defined as in
  Equation~\eqref{eqc1c2}.  Given $c_1\in C_1$ and $c_2\in C_2$, let
  \begin{equation*}
d_i(c_1)\; = \; \abs{\{x\in C_1\mid w(x-c_1) = w_i\}}\quad\text{and}\quad 
e_i(c_2)\; = \; \abs{\{x\in C_1\mid w(x-c_2) = w_i\}}
\end{equation*}
for $i \in \{1,2\}$.  Then 
\begin{equation*}
  \begin{bmatrix}
    w_1 & w_2\\ 1 & 1
  \end{bmatrix}\cdot
  \begin{bmatrix}
    d_1(c_1) \\d_2(c_1)
  \end{bmatrix} \; = \;
  \begin{bmatrix}
    D\\b_1-1
  \end{bmatrix}\qquad\text{and}\qquad
  \begin{bmatrix}
    w_1 & w_2\\ 1 & 1
  \end{bmatrix}\cdot
  \begin{bmatrix}
    e_1(c_2) \\e_2(c_2)
  \end{bmatrix} \; = \;
  \begin{bmatrix}
    E\\b_1
  \end{bmatrix},
\end{equation*}
where $D=\gamma \,n\,b_1 \,\Big[1 -
\Big(1-\frac{w_1}{\gamma\,n}\Big)^2\Big]$ and $E=\gamma \,n\,b_1
\,\Big[1 -
\Big(1-\frac{w_1}{\gamma\,n}\Big)\Big(1-\frac{w_2}{\gamma\,n}\Big)\Big]$.
In particular, the numbers $d_i(c_1),e_i(c_2)$ do not depend on the choice
of $c_1 \in C_1$ respectively $c_2\in C_2$.
\end{corollary}

Corollary \ref{corde} says that in a linear regular projective two-weight code over a
finite Frobenius ring, the number of words of weight $w_1$ at distance
$w_1$ from a word of weight $w_1$ is constant and given by $d_1 =
d_1(c_1)$, and the number of words of weight $w_1$ at distance $w_1$
from a word of weight $w_2$ is constant and given by $e_1 = e_1(c_2)$.
The linearity of $C$ then immediately gives the following result.

\begin{theorem}\label{main}
  Let $C\leq {_RR^n}$ be a regular projective two-weight code over a
  finite Frobenius ring $R$ with normalized homogeneous weights $w_1$
  and $w_2$. 
  Then the graph $\Gamma(C):=(C,E)$ with
  vertex set $C$ and edge set $E:=\bigl\{\{x,y\}\mid x, y\in
  C\text{ with }w(x-y)=w_1\bigr\}$ is strongly regular. Its
  parameters are $(N,K,\lambda,\mu)$, where
  \begin{align*}
    N & = |C|,\\
    K & = \frac{(n-w_2) |C| +w_2}{w_1-w_2}, \\
    \lambda & = \frac{n\, K
      \Big[1-\Big(1-\frac{w_1}{n}\Big)^2\Big]+w_2(1-K)}{w_1-w_2},\\
    \mu & = \frac{n\, K \Big[1-\Big(1-\frac{w_1}{
        n}\Big)\Big(1-\frac{w_2}{n}\Big)\Big]-w_2\,K}{w_1-w_2}.
  \end{align*}
\end{theorem}

We now give some simple examples illustrating the use of
Theorem~\ref{main}. Further examples will be presented in the
next section.

\begin{example}
  \label{parity-check-code}
  The code
  $C=\{x\in\Z_4^3\mid x_1+x_2+x_3=0\}$, i.e.\ the parity-check code of
  length $3$ over $\Z_4$, is a regular projective two-weight code with nonzero
  Lee weights $w_1=2$, $w_2=4$ and frequencies $b_1=6$, $b_2=9$. The
  constant weight subcodes of $C$ are
  \begin{equation*}
    C_1=\{130,310,103,301,013,031\}\quad\text{and}\quad
    C_2=\{112,332,121,323,211,233,220,202,022\}.
  \end{equation*}
  By Theorem~\ref{main}, the graph $\Gamma(C)$ obtained from the
  codewords in $C_1$ is strongly regular
  with parameters $(16,6,2,2)$. It is easily checked that the induced
  subgraph on $C_1$ (the neighborhood of the vertex $000$) is a
  $6$-cycle. Hence $\Gamma(C)$ is isomorphic to the Shrikhande graph (cf.~\cite[Ex.~4.6]{CL91}). 
  
  The other $(16,6,2,2)$ strongly regular
  graph, the square lattice graph $\mathrm{L}_2(4)$, is obtained by
  taking instead of $\Z_4$ one of the rings $\F_4$ or $\F_2[X]/(X^2)$.  
  
  Now consider fourth ring of order $4$, $\F_2\oplus\F_2$, and replace $\Z_4$ with it in the above.
  The resulting code is not proper and the constant weight subcodes containing words of weight 0 
  and 4 are
  $$C_1 = \{ ooo,occ,coc,cco\} \quad\text{and}
  \quad C_2:=\{ oaa,aoa,aao,obb,bob,bbo,abc,abc,bac,bca,cab,cba \}$$
  where $o=00$, $a = 10$, $b=01$ and $c=11$.
  Although the code is not proper, it still gives the (trivial) strongly regular with 
  parameters $(16,3,2,0)$.
  \end{example}

\section{Constructions}

The constructions given here require that the reader be familiar with the subclass of Frobenius rings called {\em chain rings}. A ring $R$ is called a {\em left chain ring\/}, if the set of its left ideals forms a chain by set inclusion. If $R$ is finite it can be seen that also the set of all right ideals of $R$ forms a chain, and hence we refer to these rings simply as finite chain rings.

The unique maximal ideal $\rad(R)$ of a finite chain ring $R$ is principal, and hence generated by an element $\theta\in R$. The nilpotency index of $\theta$ is then called the length of $R$, and our first pair of examples assume that $R$ is of length $2$. Observe that $F:=R/\rad(R)$ is a finite field, of size $q$, say, and $\rad(R)$ has $q$ elements.

On the set of points of $\barb(R^2_R)$ we define the neighbourhood relation $\sim$ by $xR\sim yR$ if and only if $\nu(xR)=\nu(yR)$ where $\nu: R^2\longrightarrow \F^2_q$ is the natural semilinear epimorphism defined by $x\mapsto x+\rad(R^2_R)$. This relation is an equivalence relation which has $q+1$ equivalence classes, and each of these classes consists of $q$ points.

\begin{proposition}\label{6.1}
  Let $R$ be a finite chain ring of length $2$ with a $q$-element residual field.
  Let $1\leq s \leq q$ be an integer. Let $Y$ be a $2 \times s(q+1)$ matrix whose   columns generate $s$ distinct elements of each equivalence class of $\sim$ given above. 
  Then $Y$ generates an $[s(q+1),2]$ two-weight code with weights $$w_1\; =\; \frac{q(qs-1)}{q-1} \quad\mbox{and}\;\; w_2\; = \; \frac{q^2s}{q-1}.$$ 
  \end{proposition}
  \pf Words of the generated code $C$ are of the form $xY$ where $x\in R^2$. If $x$ is a unimodular element (i.e.~there is $z\in R$ with $x_1z_1+x_2z_2=1$) then there are exactly $qs$ columns $z$ in $Y$ for which $xz\in R^\times$. For each of the remaining $s$ 
  columns $z$ of $Y$ we have $xz\in \rad(R)$, where there is at most one column $z$ with $xz=0$. For this reason we have $$\whom(xY) = \left\{\begin{array}{lcl}
  qs+(s-1)\frac{q}{q-1} & = & \frac{q(qs-1)}{q-1}, \;\;\mbox{or}\\
  qs+s\frac{q}{q-1} & = & \frac{q^2s}{q-1}.
  \end{array}\right.$$
  
  If $x$ is non-zero, but not unimodular, then $x=\theta x'$ for some unimodular $x'\in R^2$ (and where $\theta$ is a generator of $\rad(R)$). The codeword in question is given by $\theta x'Y$. For all columns of $Y$ generating points in the same equivalence class of $\sim$ the expression $\theta x'Y$ is constant, and it vanishes on exactly one class. This means that $w_H(xY)=qs$, and for this reason we find $\whom(xY) = qs\frac{q}{q-1}$. $C$ is now proper because $\whom$ is strictly positive. All in all we have hence proven the claim.\qed

\begin{corollary}
The code $C$ in the previous proposition determines a strongly regular graph with parameters
\[
N = q^4,\quad K=s(q^3-q),\quad \lambda = q^2(1+s^2)-3sq,\quad \mu =
sq(sq-1).
\]
\end{corollary} 

For example, let $R=\Z_9$. The ring $R$ has maximal ideal $3\Z_9$,
which has 3 elements.  Now the set of points in $\barb(\Z_9^2)$ has
$q+1=4$ distinct equivalence classes, and each class contains $3$
points. Choose $s$ points from each class, for $s = 1,2,$ or $3$ and
form a generator matrix for a $\Z_9-$code with two homogeneous weights
$w_1,w_2$. The table below gives the parameters of the codes
constructed as described before, and the parameters of the
corresponding strongly regular graph.
\begin{center}
  \begin{tabular}{|cccc|cccc|}
    \hline
    $n$ & $k$ & $w_1$ & $w_2$ & $N$ & $K$ & $\lambda$ & $\mu$ \\
    \hline
    4 & 2 & 3 & 4.5 & 81 & 24 & 9 & 6 \\
    \hline
    8 & 2 & 7.5 & 9 & 81 & 48 & 27 & 30 \\
    \hline
    12 & 2 & 12 & 13.5 & 81 & 72 & 63 & 72 \\ 
    \hline
  \end{tabular}
\end{center}
The $[12,2]$ code over $\Z_9$ with homogeneous weights 12 and 13.5 is
the simplex code, and gives a trivial strongly regular graph.

Strongly regular graphs with parameters $(81,24,9,6)$ can also be
constructed from the block graph of an orthogonal array $\OA(9,3)$, or
a projective ternary $[12,4]$ code with Hamming weights 6 and 9. The
complement of such a graph has parameters $(81,56,37,42)$ and may also
be constructed from the block graph of $\OA(9,7)$.
   
Strongly regular graphs with parameters$(81,48,27,30)$ are constructed
from the block graph of an orthogonal array $\OA(9,6)$. The
corresponding complements have parameters $(81,32,13,12)$ and may be
constructed from a projective ternary $[16,4]$ code with Hamming
weights 9 and 12 or from the block graph of an orthogonal array
$\OA(9,4)$.

The next example needs some preparation. Let $R=\F_q[X;\sigma] / (X^2)$ for some $\sigma \in\Aut(\F_q)$ be a truncated skew polynomial ring. There is a natural embedding ${\rm PGL}(3,\F_q)\stackrel{e}{\longrightarrow} {\rm PGL}(3,R)$. If $p$ is a point of $\barb(R^3_R)$ then the orbit of $p$ under an ($e$-embedded) Singer cycle of ${\rm PGL}(2,\F_q)$ is a $q^2+q+1$ element point set $K$ that shares exactly $1$ or $q+1$ points with every hyperplane of $\barb(R^3_R)$.

\begin{proposition} \label{6.2}
Let $R$ be the above truncated skew polynomial ring, and 
 let ${Y}$ be a $3 \times (q^2+q+1)$-matrix over $R$ whose columns are
 generating the points of $K$ as defined above. Then $Y$ generates a $[q^2+q+1,3]$ two-weight code with weights $$w_1\; = \; \frac{q^3}{q-1}\quad \mbox{and}\;\; w_2\; = \; q^2.$$ 
\end{proposition}
\pf We will only sketch the proof. Words of the generated code are of the form $xY$ where $x\in R^3$.
Similar to the proof of Proposition \ref{6.1} we can show that if $x$ is unimodular, then 
 $$
  \whom(xY)\; = \;\frac{q^3}{q-1} \quad \mbox{or}\;\; \whom(xY) \; = \;q^2.$$
  If $x$ is a non-zero element of $\rad(_RR^3)$ then $\whom(xY) = \frac{q^3}{q-1}$. This completes the proof.\qed
  
  \begin{corollary}
 The code $C$ in the previous proposition determines a strongly regular graph with parameters
  \begin{equation*}
  N = q^6,\quad K=q^4-q,\quad \lambda = q^3+q^2-3q,\quad \mu = q^2-q.
  \end{equation*}
\end{corollary}

  For example, for $p=2$ we get the table
  \begin{center}
    \begin{tabular}{|cccc|cccc|}
      \hline
      $n$ & $k$ & $w_1$ & $w_2$ & $N$ & $K$ & $\lambda$ & $\mu$ \\
      \hline
      7 & 3 & 4 & 8 & 64 & 14 & 6 & 2 \\
      \hline
      21 & 3 & 16 & 64/3 & 4096 & 252 & 68 & 12 \\
      \hline
      73 & 3 & 64 & 512/7 & 262144 & 4088 & 552 & 56 \\ 
      \hline
    \end{tabular}
  \end{center}
  The strongly regular graph with parameters $(64,14,6,2)$ is unique,
  and may also be constructed from a projective binary $[14,6]$ code
  with Hamming weights $4$ and $8$.

Our last example again needs further preparation. In some way it combines the ideas of the first two examples.

First recall that if $\F_q$ is the finite field of $q$ elements, then in the projective plane ${\rm PG}(F_q^3)$ there is a bijection $\sigma$ between the set of points and the set of lines, such that $p \in \sigma(p)$ for all points. $\sigma$ might for example be chosen as the orbit of a flag $(p,\ell)$ under a Singer cycle of the plane. 

Let now $R$ be a chain ring with residue field $F_q$. Similar to what we introduced earlier we have an equivalence relation $\sim$ on the set of all points of $\barb(R^3_R)$ defined by $xR\sim yR$ if and only if $\nu(xR) = \nu(yR)$. Here $\nu: R^3\longrightarrow F_q^3$ is the natural semilinear epimorphism with $x\mapsto x+\rad(R^3)_R$. This neighbourhood relation has exactly $q^2+q+1$ classes. Let $p$ be a point, and let $h$ be a hyperplane in $\barb(R^3_R)$. We define $$[p]_h \; := \; \{g\mid g \sim p\;\;\mbox{and}\;\; g | h\},$$ and call it a {\em line segment\/}, provided it is non-empty. We will call $\nu(h)$ the {\em direction\/} of the line segment.

It is now possible to compose a $q(q^2+q+1)$-element set $K$ of complemented free points of $\barb(R^3_R)$ that meets every neightbourhood class $[p]$ in a line segment in such a way that the occuring directions are all pairwise distinct. This is accomplished using the bijection $\sigma$ introduced above.

\begin{proposition}
Let $Y$ be a $3\times sq(q^2+q+1)$ matrix whose columns generate the elements of the set $K$ in $\barb(R^3_R)$. The code generated
  by $Y$ over $R$ is a two-weight code with parameters given by
  \[
  n = q(q^2+q+1),\quad k = 3,\quad w_1 = \frac{q^4-q^2}{q-1},\quad
  w_2 = \frac{q^4}{q-1},
  \]
  and determines a strongly regular graph with parameters
  \[
  N = q^6,\quad K=q^5-q^2,\quad \lambda = q^4+q^3-3q^2,\quad
  \mu = q^2(q^2-1).
  \]
\end{proposition}

\begin{remark}
This construction can be generalised to sets $K$ containing $s$ line segments of the same direction from each equivalence class of points, where $1\leq s\leq q$. Such a set $K$ will have either $wq+q$ or $sq$ points in common with every hyperplane.

The two-weight code will then have parameters given by 
 \[
  n = sq(q^2+q+1),\quad k = 3,\quad w_1 = \frac{sq^4-q^2}{q-1},\quad
  w_2 = \frac{sq^4}{q-1},
  \]
   and the induced strongly regular graph will have parameters
  \[
  N = q^6,\quad K=s(q^5-q^2),\quad \lambda = s^2q^4+q^3-3sq^2,\quad
  \mu = sq^2(sq^2-1).
  \]
  
\end{remark}
 
  For example, for $q=2$ we get codes and graphs with parameters
  \begin{center}
    \begin{tabular}{|cccc|cccc|}
      \hline
      $n$ & $k$ & $w_1$ & $w_2$ & $N$ & $K$ & $\lambda$ & $\mu$ \\
      \hline
      14 & 3 & 12 & 16 & 64 & 28 & 12 & 12 \\
      \hline
      28 & 3 & 28 & 32 & 64 & 56 & 48 & 56 \\
      \hline
    \end{tabular}
  \end{center}
  whereas for $q=3$ we get codes and graphs with parameters
  \begin{center}
    \begin{tabular}{|cccc|cccc|}
      \hline
      $n$ & $k$ & $w_1$ & $w_2$ & $N$ & $K$ & $\lambda$ & $\mu$ \\
      \hline
      39 & 3 & 36 & 40.5 & 729 & 234 & 81 & 72 \\
      \hline
      78 & 3 & 76.5 & 81 & 729 & 468 & 297 & 306 \\
      \hline
      117 & 3 & 117 & 121.5 & 729 & 702 & 675 & 702 \\ 
      \hline
    \end{tabular}
  \end{center}

\section*{Acknowledgement}
The authors feel indebted to the reviewers for a variety of remarks and corrections
which helped to improve the quality of the published paper.

\end{document}